\documentclass[twocolumn]{autart} 

\usepackage{amsfonts}
\usepackage{amsmath}
\usepackage{cancel,epsfig,graphics}
\usepackage{graphicx,subfigure}
\usepackage{mathcomp}
\usepackage{amsmath,amssymb}
\usepackage{epsfig,graphics}
\usepackage{euscript,multicol,subfigure}
\usepackage{wrapfig, rotating}
\newcommand{\bfx}{{\mathbf x}}
\newcommand{\bfA}{{\mathbf A}}
\newcommand{\bfE}{{\mathbf E}}

\newcommand{\bfI}{{\mathbf I}}
\newcommand{\bfb}{{\mathbf b}}

\newcommand{\bfc}{{\mathbf c}}

\newcommand{\bfV}{{\mathbf V}}
\newcommand{\bfW}{{\mathbf W}}

\newcommand{\res}{{\mathsf{res}}}
\newcommand{\real}{{\mathsf{Re}}}

\newcommand{\cH}{ {\mathcal H} }

\newcommand{\cHtwo}{ { {\mathcal H}_2} }
\newcommand{\cHtwoW}{ {{\mathcal H}_2\mbox{\tiny $(W)$}} }
\newcommand{\cHinf}{ { {\mathcal H}_\infty} }

\newcommand\IR{ {\mathbb R}}
\newcommand\IC{ {\mathbb C}}

\begin{document}

\begin{frontmatter}
\title{Interpolatory Weighted-$\cHtwo$ Model Reduction}

\thanks[footnoteinfo]{This paper has not been presented at any IFAC 
meeting. Corresponding author S.~Gugercin. Tel. +1-540-231-6549. 
Fax +1-540-231-5960}

\author[Anic]{Branimir Ani\'{c}}\ead{anic@kit.edu},
\author[GugBea]{Christopher Beattie}\ead{beattie@vt.edu},  
\author[GugBea]{Serkan Gugercin}\ead{gugercin@math.vt.edu},   and
\author[Antoulas]{Athanasios C. Antoulas}\ead{aca@rice.edu},    

\address[Anic]{Department of Mathematics, Karlsruhe Institute of Technology, Germany}
\address[GugBea]{Department of Mathematics, Virginia Tech, Blacksburg, VA, 24061-0123, USA}
\address[Antoulas]{Department of Electrical and Computer Engineering, Rice University, Houston, TX 77251, USA}



\begin{abstract}
This paper introduces an interpolation framework for the weighted-$\cHtwo$ model reduction problem.
We obtain a new representation of the weighted-$\cHtwo$ norm of SISO systems that provides new interpolatory first order necessary conditions for an optimal reduced-order model.   The $\cHtwo$ norm representation also provides an error expression that motivates a new weighted-$\cHtwo$ model reduction algorithm. 
Several numerical examples illustrate the
effectiveness of the proposed approach.
\end{abstract}

\end{frontmatter}

\section{Introduction}

Consider a single input/single output (SISO) linear dynamical system with a realization 
\begin{equation}
\label{ltisystemintro}
\bfE\, \dot{\bfx}(t)  =  \bfA \bfx(t) + \bfb\,u(t),~~
y(t)  = \bfc^T \bfx(t)
\end{equation}
for $\bfE,\bfA \in \IR^{n \times n}$ and $\bfb,\bfc \in \IR^{n}$.
$\bfx(t)\in \IR^n$, $u(t)\in \IR$, $y(t)\in \IR$, are respectively
the \emph{state}, \emph{input}, and \emph{output} of the system.  The \emph{transfer function} of this system is $G(s) = \bfc^T(s\bfE-\bfA)^{-1}\bfb$. 
Following common usage, the underlying system will also be denoted by $G$.
For many examples, the state-space dimension $n$ is quite large, leading to untenable demands on computational resources.
\emph{Model reduction} attempts to address this
by finding a reduced-order system  of the form,
\begin{equation} \label{redsysintro}
\bfE_r\,\dot{\bfx}_r (t)  =  \bfA_r \bfx_r (t) + \bfb_r u(t),~~
y_r(t)   =   \bfc^T_r \bfx_r (t) \quad~
\end{equation}
with $G_r(s) = \bfc_{r}^T(s\bfE_r-\bfA_{r})^{-1}\bfb_{r}$
for $\bfE_r$, $\bfA_r \in \IR^{r \times r}$ and $\bfb_r,\,\bfc_r \in
\IR^{r}$ with $r \ll n$ such that $y_{r}(t)\approx y(t)$ over a large
class of inputs $u(t)$.  $G_r$ is a low order, yet high fidelity,
approximation to $G$. 
We construct $G_r$  via  state-space projection:
Two matrices  (``reduction bases") $\bfV_r, \ \bfW_r \in \IR^{n \times r}$ are chosen.
Then, system dynamics are approximated  by $\bfx(t) \approx \bfV_r \bfx_r(t)$ and forcing a
Petrov-Galerkin conditon (``orthogonal residuals")
$
\mathbf{W}_r^{T}(\bfE \bfV_r\dot{\bfx}_r(t)-\bfA\bfV_{\! r}\bfx_r(t) -\bfb\,u(t))=\mathbf{0}$,
together with  the output equation $y_r(t) = \bfc^T\bfV_r \bfx_r(t)$
to produce
\begin{equation}  \label{eqn:reduction}
\begin{array}{c}
\bfE_r = \bfW_r^T \bfE \bfV_r, \quad
\bfb_r = \bfW_r^T\bfb,  \\
\bfA_r = \bfW_r^T \bfA \bfV_r, \quad \mbox{and} \quad
\bfc_r^T = \bfc^T\bfV_r.
\end{array}
\end{equation}
See \cite{An05,antoulas2010imr} for more information.
\subsection{Model Reduction by Interpolation}
The reduction bases, $\bfV_r$ and $\bfW_r$, used in (\ref{eqn:reduction}) will be chosen to force interpolation:  $G_r(s)$ will interpolate $G(s)$ (possibly together with higher order derivatives) at selected interpolation points. 
This approach to rational interpolation has been considered in \cite{yousouff1984cer,yousuff1985lsa,devillemagne1987mru,grimme1997krylov,feldmann1995efficient,antoulas2010imr} 
and depends on the following result.
\begin{thm} \label{interpCor}
Given two sets of interpolation points $\{\sigma_k\}_{k=1}^r$ and $\{\zeta_k\}_{k=1}^r$, that are each closed under conjugation, and a dynamical system $G$ as in (\ref{ltisystemintro}),   consider matrices  $\bfV_{\! r}$ and $\mathbf{W}_{\! r}$ such that
\vspace{-4ex}
{\small
\begin{equation} \label{interpCond}
\begin{array}{l}
\mbox{Range}(\bfV_{\! r}) = {span}
\left\{(\sigma_i\bfE-\bfA)^{-1}\bfb\right\}\
 \\
\mbox{Range}(\mathbf{W}_{\! r}) = {span}
\left\{({\zeta_i}\bfE-\bfA)^{-T}\bfc\right\}
\end{array}
~~\mbox{for}~~ i=1,\ldots,r.
\end{equation}
}
Then, $\bfV_{\! r}$ and $\mathbf{W}_{\! r}$ can be chosen to be real;  $G_r(s) =  \bfc_r^T(s\bfE_r- \bfA_r)^{-1}\bfb_r$ defined by (\ref{redsysintro})-(\ref{eqn:reduction})
is a real dynamical system that satisfies
$ G(\sigma_k) = G_r(\sigma_k)$ and $G(\zeta_k) = G_r(\zeta_k)$ for $k=1,\ldots,r$;
and, if $\sigma_j=\zeta_j$ for some $j$, then
$
G'(\sigma_j) = G_r'(\sigma_j),
$ 
as well where $G'$ denotes the derivative of $G(s)$ with respect to $s$.
\end{thm}

Theorem \ref{interpCor}  can be  generalized to  higher-order derivative interpolation as well, see \cite{yousouff1984cer,yousuff1985lsa,devillemagne1987mru,grimme1997krylov,feldmann1995efficient,antoulas2010imr}.  
The subspaces 
of Theorem \ref{interpCor} are  \emph{rational Krylov subspaces} and so, interpolatory model reduction methods for SISO systems are sometimes referred  to as  \emph{rational Krylov methods}. 
\subsection{Weighted Model Reduction}
The $\cHinf$ norm of a stable linear system associated with a transfer function, 
$G(s)$, is defined as 
$\displaystyle \left\| G \right\|_{\cHinf}=
\max_{\omega\in {\mathbb R}}\left|G(\imath \omega)\right|.$
The $\cHtwo$ norm of $G$ is defined as 
$\displaystyle \left\| G \right\|_{\cHtwo} :=
\left(\frac{1}{2\pi}\int_{-\infty}^{\infty}
\left| G(\imath \omega )  \right|^2 d\omega\right)^{1/2}.$ The vector spaces of meromorphic functions that are analytic in the right halfplane, having either bounded $\cHinf$ norm or bounded $\cHtwo$ norm will be denoted simply as $\cHinf$ or $\cHtwo$, respectively.   Let $W\in\cHinf $ be given.  The ($W$-)\emph{weighted $\cHtwo$ norm} is defined as $ \| G  \|_{\cHtwoW}= \| G \cdot W \|_{\cHtwo}. $

We are interested  in finding a reduced-order model $G_r$ that minimizes a $W$-weighted $\cHtwo$ norm, i.e., that solves
\begin{equation} \label{H2optProb}
 \| G - G_r \|_{\cHtwoW} =  
 \min_{ \dim(\widetilde{G}_{r})=r } \| G - \widetilde{G}_{r} \|_{\cHtwoW}
\end{equation}
The introduction of $W(s)$ allows one to penalize the error in certain frequency ranges more heavily than in others. 


{\bf An illustrative example: controller reduction}
Consider a linear dynamical system,  $P$ (the \emph{plant}) together with an associated stabilizing controller,  $G$,  
that is connected to $P$ in a feedback loop.
Many control design  methodologies,
such as LQG and $\cHinf$ methods, lead ultimately to
controllers whose order is generically as high as the order of the plant,
see \cite{varga2003aem,ZhDG96} and references therein.  Thus, high-order plants generally
lead to high-order controllers.  However,
high-order controllers are usually undesirable in real-time applications due to  complex hardware,  
 degraded accuracy,  and degraded computational speed. Thus, one prefers to use a reduced controller
$G_r$  to replace $G$.
Requiring $G_r$ to be a good approximation to $G$ is often not enough
in terms of closed-loop performance; plant dynamics need to be taken into account during the reduction process. This may be achieved through frequency
weighting:  Given a  stabilizing controller $G$, if 
$G$ has the same number of unstable poles as $G_r$ and  if
$\left\| [G-G_r]P[I + PG]^{-1} \right\|_\cHinf < 1$, then $G_r$ will also be a stabilizing controller  \cite{anderson2002crc,ZhDG96}.
Hence the controller reduction problem may be formulated as finding a reduced-order controller $G_r$ that minimizes or reduces  the weighted error $\left\| (G-G_r) W \right\|_\cHinf$ with $W(s):=P(s)(I+P(s)G(s))^{-1}$; i.e., controller reduction becomes an application of  weighted model
reduction.   This approach has been considered in 
\cite{varga2003aem,anderson2002crc,schelfhout2002ncl,gugercin2004smr,enns1984mrw,wang2002bpp,lin1992mrv,wang1999nfw,sreeram2005fwm} and references therein, leading  to variants of frequency-weighted balanced truncation. 
Conversely, the methods in \cite{halevi1992fwm} and  \cite{spanos1992anewalgorithm} 
are tailored instead towards minimizing a weighted-$\mathcal{H}_2$ error as in (\ref{H2optProb}). 

\section{Weighted-$\cHtwo$ model reduction} \label{sec:weightedh2}

The  methods proposed in  \cite{halevi1992fwm} and  \cite{spanos1992anewalgorithm} for approaching (\ref{H2optProb}) require solving a sequence of large-scale Lyapunov or Riccati equations; they rapidly become computationally intractable as the system order, $n$, increases. 
We will approach (\ref{H2optProb}) within an interpolatory model reduction framework requiring only the solution of (generally sparse) linear systems and no  need for dense matrix computations or solution of large-scale Lyapunov or Riccati equations.  Interpolatory approaches can be effectively applied  even when 
$n$
 reaches the tens of thousands.  

\subsection{A representation of the weighted-$\cHtwo$ norm}

Given transfer functions $G$,  $H\in \cHtwo$, and $W\in\cHinf$, 
define the weighted-$\cHtwo$ inner product as
\vspace{-4ex}
{\small
\begin{align*} \label{h2innerprod}
\left\langle G,\ H \right\rangle_{\cHtwoW}
&= \frac{1}{2\pi} \int_{-\infty}^{\infty}
 \overline{G(\imath \omega)W(\imath\omega)} W(\imath\omega)H(\imath \omega) \,d\omega \\
 & = \frac{1}{2\pi} \int_{-\infty}^{\infty}
 G(-\imath \omega)W(-\imath\omega)W(\imath\omega) H(\imath \omega) \,d\omega,
\end{align*}
}

\vspace{-3ex}
\noindent
 so that 
 $\displaystyle  \| G \|_{\cHtwoW} = \sqrt{\left\langle G,\ G \right\rangle_{\cHtwoW}}$.
The following lemma gives a compact expression  for 
the weighted-$\cHtwo$ inner product based on the poles and residues of 
$G(s)$, $H(s)$ and $W(s)$. By $\res[M(s),\pi]$, we denote the residue of $M(s)$
at  $\pi \in \IC$.
\begin{lem}\label{innerprodresweight}
Suppose  $G$,  $H\in \cHtwo$  have poles denoted respectively as 
$\{\lambda_{1},\dots,\lambda_{n}\}$ and $\{\mu_{1},\dots,\mu_{m}\}$,
and suppose $W\in\cHinf$ has poles denoted as 
 $\{\gamma_{1},\dots,\gamma_{p}\}$.  Assume that 
   $H(s)$ and $W(s)$ have no common poles, and 
 the poles of $W(s)$ are simple. Then  
\begin{align*} 
\left\langle G,\ H \right\rangle_{\cHtwoW} &= 
\sum^{m}_{k=1} \res[G(-s)W(-s)W(s)H(s),\mu_{k}] \\ &~~+\sum^{p}_{i=1}
G(-\gamma_{i})W(-\gamma_{i})H(\gamma_{i})\cdot\res[W(s),\gamma_{i}].
\end{align*}
Define $\chi_k =\res[G(-s)W(-s)W(s)H(s),\mu_{k}] $.
\begin{itemize}
\item If $\mu_{k}$ is a simple pole of $H(s)$, then
\begin{align*}
\chi_k{=}G(-\mu_{k})W(-\mu_{k})W(\mu_{k})\cdot \res[H(s),\mu_{k}].
\end{align*}
\item If $\mu_{k}$ is a double pole of $H(s)$, then
\begin{align*}
\chi_k= & G(-\mu_{k})W(-\mu_{k})W(\mu_{k})\cdot\res[H(s),\mu_{k}]\\ 
&-h_{-2}(\mu_{k})\,\cdot\,\frac{d}{ds}\left.\left[G(s)W(s)W(-s)\right]\right|_{s=-\mu_{k}},
\end{align*}
where
$\displaystyle h_{-2}(\mu_{k}) = \lim_{s\to\mu_{k}}(s-\mu_{k})^{2}H(s)$.
\end{itemize}
\end{lem}

 {\bf Proof:}
 $T(s) = G(-s)W(-s)W(s)H(s)$ has poles at
 $$
 \{-\lambda_1,\ldots\,-\lambda_n\} \cup \{\pm\gamma_1,\ldots,\pm\gamma_p\} \cup \{\mu_1,\dots,\mu_m\}.
 $$
For any  $R>0$,   define a semicircular contour in the left
halfplane:
$
\Gamma_{R}=\left\{z\left| z=\imath\omega \mbox{ with
}\omega\in[-R,R]\right. \right\}\cup \left\{z\left| z=R\,
e^{\imath\theta} \mbox{ with
}\theta\in[\frac{\pi}{2},\frac{3\pi}{2}]\right. \right\}.
$
For $R$ large enough, the region bounded by $\Gamma_{R}$ contains 
$
\{\gamma_1,\ldots,\gamma_p\}\cup\{\mu_1,\dots,\mu_m\},
$
 constituting all the poles of $W(s)H(s)$, and hence all the stable poles of $G(-s)W(-s)W(s)H(s)$.
 Then, the Residue Theorem  yields

\vspace{-4ex}
{\small
 \begin{align*}
\left\langle G,\ H \right\rangle_{\cHtwoW} &=
\frac{1}{2\pi} \int_{-\infty}^{+\infty}
G(-\imath \omega)W(-\imath \omega)W(\imath \omega) H(\imath \omega) \,d\omega \\
& =  \lim_{R\rightarrow\infty} \frac{1}{2\pi\imath} \int_{\Gamma_{R}}G(-s)W(-s)W(s) H(s)\  ds \\
  &=   \sum_{k=1}^{m} \mbox{\textsf{res}}[G(-s)W(-s)W(s) H(s),\mu_k] \\
&\qquad+\sum_{i=1}^{p} \mbox{\textsf{res}}[G(-s)W(-s)W(s) H(s),\gamma_i].
\end{align*}
}

\vspace{-4ex}
\noindent
This leads to the first assertion.   
If $\mu_{k}$ is a simple pole for $H(s)$, then  
\begin{align*}
& \mbox{\textsf{res}}[G(-s)W(-s)W(s)H(s),\,\mu_{k}] \\
& =
\lim_{s\rightarrow \mu_{k}} [(s-\mu_{k})G(-s)W(-s)W(s) H(s)]\\
 &= G(-\mu_{k})W(-\mu_{k})W(\mu_{k})\lim_{s\rightarrow \mu_{k}} (s-\mu_{k})H(s).
\end{align*}
Similarly, if $\mu_{k}$ is a double pole for $H(s)$, then it is also a double
pole for $G(-s)W(-s)W(s)H(s)$ and
 \begin{align*}
\mbox{\textsf{res}}[&G(-s)  W(-s)W(s) H(s),\,\mu_{k}] \\
&= \lim_{s\rightarrow \mu_{k}} \frac{d}{d s} [(s-\mu_{k})^{2}G(-s)W(-s)W(s) H(s)] \\
 &= \lim_{s\rightarrow \mu_{k}} G(-s) W(-s)W(s)\frac{d}{d s}\left[(s-\mu_{k})^{2} H(s)\right]\\ 
& ~~~+ \lim_{s\rightarrow \mu_{k}}(s-\mu_{k})^{2} H(s)\frac{d}{d s}\left[G(-s)W(-s)W(s)\right]\\
&=G(-\mu_{k})W(-\mu_{k}) W(\mu_{k})\,\cdot\,\mbox{\textsf{res}}[H(s),\,\mu_{k}]\\
&~~~-h_{-2}(\mu_{k})\,\cdot\,\frac{d}{d s}\left.\left[G(s)W(s)W(-s)\right]\right|_{s=-\mu_{k}}
~\Box
\end{align*} 
\vspace{-6ex}
\begin{cor}\label{normHtwo}
If $G(s)$ and $W(s)$ in Lemma \ref{innerprodresweight} each have simple poles, then 
\vspace{-6ex}
{ 
\begin{eqnarray}  \| G \|_{\cHtwoW}^{2}  = \sum_{k=1}^{n}
G(-\lambda_{k})W(-\lambda_{k})W(\lambda_{k})\cdot\mbox{\textsf{res}}[G(s),\lambda_k]  \nonumber \\
+\sum_{k=1}^{p}G(-\gamma_{k})W(-\gamma_{k})G(\gamma_{k})\cdot\mbox{\textsf{res}}[W(s),\gamma_k]. \qquad
 \label{eq:wh2norm}
\end{eqnarray}
}
\end{cor}

\vspace{-4ex}
This new formula (\ref{eq:wh2norm}) for the weighted-$\cHtwo$ norm contains as a special case (with $W(s) = 1$), 
a similar expression for the (unweighted) $\cHtwo$ norm introduced in \cite{gugercin2008hmr}.

Suppose $W\in\cHinf$ has simple poles at 
 $\{\gamma_{1},\dots,\gamma_{p}\}$ and 
 define a linear mapping $\mathfrak{F}:\cHtwo\rightarrow \cHtwo$ by
 \vspace{-6ex}
 {\small
\begin{equation} \label{Fmap}
\mathfrak{F}[G](s)=  G(s)W(s)W(-s)+\sum^{p}_{k=1}
G(-\gamma_{k})W(-\gamma_{k})\frac{\res[W(s),\gamma_{k}]}{s+\gamma_{k}}
\end{equation}
}
Notice that $G(s)W(s)W(-s)$ has simple poles at $-\gamma_{1},-\gamma_{2},\dots,-\gamma_{p},$ and 
\begin{align*}
\res[G(s)W(s) & W(-s),-\gamma_{k}]\\ & =\lim_{s\rightarrow -\gamma_{k} }(s+\gamma_{k})G(s)W(s)W(-s) \\
 &= G(-\gamma_{k})W(-\gamma_{k})\lim_{s\rightarrow -\gamma_{k} }(s+\gamma_{k})W(-s) \\
&  = -G(-\gamma_{k})W(-\gamma_{k})\lim_{s\rightarrow \gamma_{k} }(s-\gamma_{k})W(s) \\
 & =-G(-\gamma_{k})W(-\gamma_{k})\cdot\res[W(s),\gamma_{k}].
\end{align*}
Thus $\mathfrak{F}[G](s)$ has poles only in the left half plane and indeed $\mathfrak{F}:\cHtwo\rightarrow \cHtwo$.
\begin{cor}\label{ElemOrthCond}
	Suppose  $G$ and $W$ are stable with  poles
$\{\lambda_{1},\dots,\lambda_{n}\}$ and 
 $\{\gamma_{1},\dots,\gamma_{p}\}$, respectively. Choose $\mu$ arbitrarily in the left half plane distinct from these points.
 Then for $F(s)=\mathfrak{F}[G](s)$,
$ \left\langle G,\ \frac{1}{s-\mu} \right\rangle_{\cHtwoW} = F(-\mu)$  and $\left\langle G,\ \frac{1}{(s-\mu)^2} \right\rangle_{\cHtwoW} = -F'(-\mu)$.
\end{cor}
{\bf Proof:}
By Lemma \ref{innerprodresweight}, 
\vspace{-5ex}
{\small 
\begin{eqnarray*}
 \left\langle G,\ \frac{1}{s-\mu} \right\rangle_{\cHtwoW} = G(-\mu)W(-\mu)W(\mu) \qquad\qquad\qquad\\ \qquad\qquad+\sum^{p}_{k=1}
G(-\gamma_{k})W(-\gamma_{k})\frac{\res[W(s),\gamma_{k}]}{\gamma_{k}-\mu} = F(-\mu),~ {\rm and}
\end{eqnarray*}
}
\vspace{-7ex}
{\small
\begin{eqnarray*}
 \left\langle G,\ \frac{1}{(s-\mu)^2} \right\rangle_{\cHtwoW} = 
-\frac{d}{ds}\left.\left[G(s)W(s)W(-s)\right]\right|_{s=-\mu} \\
 +\sum^{p}_{k=1}
G(-\gamma_{k})W(-\gamma_{k})\frac{\res[W(s),\gamma_{k}]}{(\gamma_{k}-\mu)^2} = -F'(-\mu).~\Box
\end{eqnarray*}
}
\vspace{-2ex}
\subsection{Weighted-$\cHtwo$  optimality conditions}    \label{sec:weightedh2cond}

Consider the problem of finding a  reduced order system, $G_r$, that solves (\ref{H2optProb}).
The feasible set 
for (\ref{H2optProb}) is nonconvex, so finding a true (global) minimizer is generally intractable.  Nonetheless, we are able to obtain descriptive \emph{necessary conditions} for $G_r$ to satisfy (\ref{H2optProb}). 

\begin{thm}
If $G_r$ has simple poles, $\{\hat{\lambda}_1,\,\ldots,\,\hat{\lambda}_r\}$, and solves (\ref{H2optProb}), then $G_r$ must satisfy: for $k=1,\,\ldots,\,r$,
\begin{equation}\label{orthCondweight}
F_r(-\hat{\lambda}_k)=F(-\hat{\lambda}_k) \quad\mbox{and}\quad F_r'(-\hat{\lambda}_k)=F'(-\hat{\lambda}_k)
\end{equation}
where $F=\mathfrak{F}[G]$ and $F_r=\mathfrak{F}[G_r]$ are defined from (\ref{Fmap}).
\end{thm}
{\bf Proof:}
Suppose by way of contradiction that, for some $\mu\in \{\hat{\lambda}_1,\,\ldots,\,\hat{\lambda}_r\}$, 
$
\left\langle G-G_r,\ \frac{1}{s-\mu} \right\rangle_{\!\cHtwoW} \!\!= \alpha_0 \neq 0.
$
By hypothesis, $G_r$ can be represented as $G_r(s)=\sum_{i=1}^r \frac{\hat{\varphi}_i}{s-\hat{\lambda}_i}$ and for some index $k$,  
$\mu=\hat{\lambda}_k$.   Define $\vartheta_0=\arg(\alpha_0)$ and with $\varepsilon>0$, define
$$
\widetilde{G}_{r}^{(\varepsilon)}(s)= \frac{\hat{\varphi}_k+\varepsilon\,e^{-\imath \vartheta_0}}{s-\mu} 
+ \sum_{i\neq k} \frac{\hat{\varphi}_i}{s-\hat{\lambda}_i}.
$$
Then 
\vspace{-5ex}
{\small $$ \|G_r-\widetilde{G}_{r}^{(\varepsilon)}\|_{\cHtwoW}=\left\|\frac{-\varepsilon\,e^{-\imath \vartheta_0}}{s-\mu}\right\|_{\cHtwoW} 
\leq \|W\|_{\cHinf}\frac{\varepsilon}{\sqrt{2|\real(\mu)|}}$$
}
 so that 
$\| G_r(s)-\widetilde{G}_{r}^{(\varepsilon)}(s) \|_{\cHtwoW} =\mathcal{O}(\varepsilon)$ as $\varepsilon\rightarrow 0$.
Since $G_r$ solves (\ref{H2optProb}), 
\begin{align*}
    \| G - &G_r \|_{\cHtwoW}^{2}\leq
    \| G - \widetilde{G}_{r}^{(\varepsilon)} \|_{\cHtwoW}^{2}  \\
    \leq &\| (G - G_r) + (G_r- \widetilde{G}_{r}^{(\varepsilon)})  \|_{\cHtwoW}^{2}\\
    \leq&\| G - G_r \|_{\cHtwoW}^{2}
    + 2\real \left\langle G-G_r,\ G_r-\widetilde{G}_{r}^{(\varepsilon)}\right\rangle_{\!\cHtwoW}
 \\ &~~~~+\| G_r- \widetilde{G}_{r}^{(\varepsilon)} \|_{\cHtwoW}^{2}.
\end{align*}
Thus, \\$\displaystyle
0\leq 2\ \real \left\langle G-G_r,\
G_r-\widetilde{G}_{r}^{(\varepsilon)}\right\rangle_{\!\cHtwoW} +
\| G_r - \widetilde{G}_{r}^{(\varepsilon)} \|_{\cHtwoW}^{2}. $
This implies first that $0\leq  -\varepsilon |\alpha_0| +\mathcal{O}(\varepsilon^2)$, which then leads to a contradiction, $\alpha_0=0$. 

To show the next assertion, suppose that for some $\mu\in \{\hat{\lambda}_1,\,\ldots,\,\hat{\lambda}_r\}$, 
$
 \left\langle G-G_r,\ \frac{1}{(s-\mu)^2} \right\rangle_{\!\cHtwoW} \!\!= \alpha_1 \neq 0.
$
Then for some $k$,  $\mu=\hat{\lambda}_k$ and we define $\vartheta_1=\arg(\hat{\varphi}_k\cdot \alpha_1)$. 
For $\varepsilon>0$  sufficiently small, define
$$
\widetilde{G}_{r}^{(\varepsilon)}(s)= \frac{\hat{\varphi}_k}{s-(\mu+\varepsilon\,e^{-\imath \vartheta_1})} 
+ \sum_{i\neq k} \frac{\hat{\varphi}_i}{s-\hat{\lambda}_i} 
$$
As $\varepsilon\rightarrow 0$, we have
$$ \|G_r-\widetilde{G}_{r}^{(\varepsilon)}\|_{\cHtwoW}=\left\|\frac{-\varepsilon\,\hat{\varphi}_k\,e^{-\imath \vartheta_1}}{(s-\mu)^2-\varepsilon\,e^{-\imath \vartheta_1}}\right\|_{\cHtwoW} 
=\mathcal{O}(\varepsilon)
$$
Following a similar argument as before, we find that $0\leq  -\varepsilon |\hat{\varphi}_k\cdot \alpha_1| +\mathcal{O}(\varepsilon^2)$ as 
$\varepsilon\rightarrow 0$, 
which leads to a contradiction, $\alpha_1=0$.
$\Box$

The interpolation conditions described in (\ref{orthCondweight}) give first order necessary conditions for $G_r$ to solve the optimal weighted-$\cHtwo$ model reduction problem (\ref{H2optProb}).  Note that for $W(s) =1$,  one obtains $F(s) = G(s)$ and $F_r(s) = G_r(s)$; thus (\ref{orthCondweight}) contains the interpolatory $\cHtwo$ optimality conditions of \cite{gugercin2008hmr} for the unweighted problem as a special case.
Unfortunately, there does not appear to be a straightforward generalization of the corresponding computational approach that was described in \cite{gugercin2008hmr} for the optimal (unweighted) $\cHtwo$ model reduction problem.   Instead, we consider a different systematic approach to this problem motivated by an expression for the weighted-$\cHtwo$ error.

\subsection{A weighted-$\cHtwo$ error expression}

The weighted-$\cHtwo$ norm expression  in Corollary \ref{normHtwo}  leads  immediately to 
an expression for the weighted-$\cHtwo$ error that  forms the basis for our computational approach. 
\begin{cor}  \label{thm:h2error}
Suppose that $G$, 
$G_r$ and  $W$ are stable with simple poles $\left\{ \lambda_i\right\}_{i=1}^n$, 
$\left\{\hat{\lambda}_{j}\right\}_{j=1}^r,$
and $\left\{ \gamma_k\right\}_{k=1}^p$, respectively, and that there are no common poles.
 Define residues: 
$\phi_{i}$:= $\mbox{\textsf{res}}[G(s),\lambda_i]$; 
$\hat{\phi}_{j}$:= $\mbox{\textsf{res}}[G_{r}(s),\hat{\lambda}_j]$; and 
$\psi_{k}$:= $\mbox{\textsf{res}}[W(s),\gamma_k]$.
 The weighted-$\cHtwo$
error is given by
\vspace{-5ex}
{\small
\begin{align}
 \nonumber \big\| G&- G_{r} \big\|^{2}_{\cHtwoW} = \sum_{i=1}^{n}(G(-\lambda_{i})-G_{r}(-\lambda_{i}))W(-\lambda_{i})W(\lambda_{i})\cdot\phi_{i}\\
\label{eq:h2error} +& \sum_{j=1}^{r}(G_{r}(-\hat{\lambda}_{j})-G(-\hat{\lambda}_{j}))W(-\hat{\lambda}_{j})W(\hat{\lambda}_{j})\cdot\hat{\phi}_{j}\\
\nonumber &+ \sum_{k=1}^{p}(G(-\gamma_{k})-G_{r}(-\gamma_{k}))W(-\gamma_{k})(G(\gamma_{k})-G_{r}(\gamma_{k}))\cdot\psi_{k}
\end{align}
}
\end{cor}

\vspace{-2ex}
One may recover the (unweighted) $\cHtwo$ error expression of \cite{gugercin2008hmr} as a special case by taking $W(s)  = 1$.
Notice that the weighted error depends on the mismatch of $G$ and $G_r$ at the reflected full system poles 
$\{-\lambda_i\}$, reflected reduced poles $\{-\hat{\lambda}_j\}$,  and reflected weight poles $\{-\gamma_k\}$. 

\vspace{-1ex}
 \subsection{An algorithm for the weighted-$\cHtwo$ model reduction problem: W-IRKA}
 \vspace{-2ex}
 
 In order to reduce the weighted error, one may eliminate some terms in the error expression, by forcing {\it interpolation} at selected (mirrored) poles.
Since $r$
is required to be much smaller than $n$, there is not enough degrees of freedom to force interpolation at all the terms in the first and third components of the weighted-$\cHtwo$ error. However, the second term, i.e.
the mismatch  at $\hat{\lambda}_j$,  can be completely
eliminated by enforcing $G(-\hat{\lambda}_j) = G_r(-\hat{\lambda}_j)$
for $j=1,\ldots,r$. Hence, as in the unweighted $\cHtwo$ problem, the mirror images of the reduced-order poles play a crucial role. This motivates an algorithm with iterative rational Krylov steps to enforce the desired interpolation property as outlined in Algorithm 1 below. However, a crucial difference from the unweighted $\cHtwo$ problem is that we will not enforce interpolation of $G'(s)$ at these points; instead we use the remaining $r$ degrees of freedom to reflect the weight information $W(s)$ and also to eliminate terms from the first component of the error term.  The error expression (\ref{eq:h2error}) shows that interpolation errors are multiplied by the residues $\phi_{i}$ and $\psi_{k}$. Hence, we use the remaining $r$ variables to eliminate terms in the first and third components of the error expression corresponding to the 
dominant residues $\phi_{k}$ and $\psi_{k}$. 
Note that in several cases, such as in the controller reduction problem, the 
state-space dimension of the weight will be of the same order as that of  $G$, $\mathcal{O}(p) \approx \mathcal{O}(n)$. 
We measure dominance in a relative sense; i.e., normalized by the largest (in amplitude) $\phi_{k}$ and $\psi_{k}$ in every set. More details on this selection process can be found in Section \ref{sec:numerics} where several examples are used to illustrate these concepts. 
Note that one {\it never} needs to compute a full eigenvalue
decomposition to obtain the residues of $G(s)$ and $W(s)$. Since only a small subset of poles is needed, one could use, for example, the
{\it dominant pole algorithm} proposed by Rommes \cite{rommart06} which  computes  effectively those eigenvalues   that correspond to the dominant residues  without requiring a full eigenvalue decomposition. 
\begin{figure}
\framebox[3.35in][t]{
  \hspace*{-.1in}\begin{minipage}[l]{3.25in}
  {
 {\small
 \textbf{Algorithm 1. Weighted Iterative Rational Krylov Algorithm (W-IRKA)}
 
 \vspace{-2ex}
\begin{enumerate}
\item[] Given $G(s) = \bfc^T(s \bfE- A\bf)^{-1}\bfb$ and $W(s) = \bfc_w^T(s \bfE_w - \bfA_w)^{-1}\bfb_w$, reduction order $r=\nu+\varpi$ with $\nu,\varpi\geq 0$, let  $\{ \lambda_i\}_{i=1}^\nu$ denote  the $\nu$ dominant poles of $G$  and  $\{ \gamma_k\}_{k=1}^\varpi$ the $\varpi$
dominant poles of $W$.
  \vspace{-1ex}
 \item  Make an initial interpolation point selection:\\
 $
 \zeta_i = -\lambda_i\ \mbox{ for}~i=1,\ldots,\nu, \quad \zeta_{j+\nu} = -\gamma_{j}\ \mbox{ for } j=1,\ldots,\varpi; \quad \sigma_k=\zeta_k\ \mbox{ for } k=1,\ldots,r;
 $
  \vspace{-1ex}
\item Construct bases, $\bfV_{\! r}$ and $\mathbf{W}_{\! r}$,  that satisfy (\ref{interpCond}).
  \vspace{-1ex}
\item Repeat, while (relative change in $\{\sigma_i\} > \mbox{tol}$)
\begin{enumerate}
  \vspace{-1ex}
\item $\bfA_r = \mathbf{W}_{\! r}^T \bfA \bfV_{\! r}$ and $\bfE_r = \mathbf{W}_{\! r}^T \bfE \bfV_{\! r}$
\item Solve the eigenvalue problem $\bfA_r \bfx_j = \hat{\lambda}_j \bfE_r \bfx_j$ and assign 
$\sigma_j \longleftarrow -\hat{\lambda}_j$ for $j=1,\ldots,r$.
\item Update $\bfV_{\! r}$ so that
{$\mbox{Range}(\bfV_{\! r}) = \mbox{span}\left\{(\sigma_1 \bfE - \bfA)^{-1}\bfb,\ \cdots,~(\sigma_r \bfE - \bfA)^{-1}\bfb\right\}$}.
\end{enumerate}
  \vspace{-1ex}
\item $\bfA_r = \mathbf{W}_{\! r}^T \bfA \bfV_{\! r}$, $\bfE_r = \mathbf{W}_{\! r}^T \bfE \bfV_{\! r}$, $\bfb_r = \mathbf{W}_{\! r}^T \bfb$, and
$\bfc^T_r = \bfc^T\bfV_{\! r}$
\end{enumerate}
}
  }
 \end{minipage}
  }
\end{figure}

 Upon convergence of Algorithm 1, $\sigma_j = -\hat{\lambda}_j$ for $j=1,\ldots,r$; $G_r$  interpolates $G$ at these points, and the second sum in  (\ref{eq:h2error}) is eliminated.  $\bfW_r$ is unchanged throughout, so $G_r$ interpolates $G$ at $r$ (aggregated) dominant poles of $G$ and $W$, eliminating $\nu$ and $\varpi$ terms from the first and third sums in (\ref{eq:h2error}), respectively.  Examples in Section \ref{sec:numerics} illustrate the effectiveness of this approach.

\vspace{-4ex}
\section{Numerical examples}       \label{sec:numerics}
\vspace{-2ex}
We provide two  examples  related to controller reduction. 
 $\Phi^{(N)}$ and $\Psi^{(N)}$ denote the set of normalized
residues
of $G(s)$ and $W(s)$, respectively. 

\subsection{A building model}
\vspace{-2ex}
The plant, $P$, is linearized a model for the Los Angeles University
Hospital, and has order $48$; see \cite{antoulas2001smr} for details. 
An LQG-based controller, $G$, of the same order, $n=48$, is designed to dampen oscillations in the impulse response.
 The ten highest normalized residues of $G(s)$ and of $W(s)$ are:
\vspace{-2ex}
\begin{align*}
\Phi^{(N)} =  [~&1.0000 ~1.0000 ~0.0286 ~0.0286~
0.0088, \\
& 0.0088~0.0080 ~0.0080 ~0.0060 ~0.0060~] \\
\Psi^{(N)} = [~&1.0000 ~1.0000 ~0.8416 ~0.8416 ~0.3935,\\
&0.3935~0.2646~0.2646 ~0.0951 ~0.0951~]
\end{align*}

\vspace{-4ex}
There is a significant drop in $\Phi^{(N)}$ values after the second entry, so we take the first two residues of $G$ as dominant.  $\Psi^{(N)}$ remains
at roughly the same order until the $9^{\rm th}$ entry.
Thus, we choose $\nu = 2$; and  
$\varpi = r - \nu = r-2$ for a given reduction order, $r$. 
To illustrate the effect of this dominant pole selection, 
we apply {\bf W-IRKA}, varying $\nu$ from $0$ to $r$. 
Tables \ref{table-48-16} below lists the resulting weighted-$\cHtwo$ errors 
for three cases: $r=12$, $r=14$, and $r=16$.
\begin{table}[hh]
\centering 
$r=12$:\\ 
\vspace{1ex}
$
 {\scriptsize
  \begin{array}{|c|c|c|c|c|c|c|c|c|}
    \hline
    \nonumber\mbox{$\nu/ \varpi$} & \mbox{$12/0$} & \mbox{$10/2$} & \mbox{$8/4$} & \mbox{$6/6$} & \mbox{$4/8$}& \mbox{$2/10$} & \mbox{$0/12$}  \\\hline
  \multicolumn{1}{c|}  ~ & 1.4021 & 1.1433 & 0.6548 & 0.6863 & 0.3576 & 0.2181 & 0.2853 \\
  \cline{2-8}
     \end{array}
     }
\vspace{2ex}
$
 $r=14:$\\  
 \vspace{1ex}
$
  {\scriptsize\begin{array}{|c|c|c|c|c|c|c|c|c|c|}
    \hline
    \nonumber
    \mbox{$\nu/ \varpi$} & \mbox{$14/0$} & \mbox{$12/2$} & \mbox{$10/4$} & \mbox{$8/6$} & \mbox{$6/8$} & \mbox{$4/10$} & \mbox{$2/12$} & \mbox{$0/14$} \\\hline
    \multicolumn{1}{c|}  ~ & 1.4734 &1.3436 & 0.6477 & 0.3019 & 0.1538 & 0.1425 & 0.1351 & 0.2224  \\
  \cline{2-9}
  \end{array}
}  \vspace{2ex}
$
$r=16:~~~$   
$
  {\scriptsize 
  \begin{array}{|c|c|c|c|c|c|}
 \hline
    \nonumber\mbox{$\nu/\varpi$} & $16/0$ & $14/2$& $12/4$
    & $10/6$ & $8/8$ 
     \\\hline
     \multicolumn{1}{c|}  ~ & 1.4206 & 1.1934 &  0.7258 & 0.2898
  & 0.1917      \\\hline
      \nonumber\mbox{$\nu/\varpi$} &  $6/10$& $4/12$ & $2/14$& $0/16$
     \\\cline{1-5}
    \multicolumn{1}{c|}  ~  & 0.1221& 0.1154 &  0.1309 & 0.1388
     \\\cline{2-5}

    \end{array}
    }
$ 
\caption{Weighted-$\cHtwo$ error as $\nu$ and $\varpi$ vary} 
\label{table-48-16}
\end{table}

The weighted-$\cHtwo$
error decreases as we take more dominant poles of $W(s)$ over those of 
$G(s)$; suggesting the importance of the residues of $W(s)$ in the error expression (\ref{eq:h2error}).  
Choosing $\nu$=2 is the best choice for
most cases.  Tables \ref{table-48-16}  illustrate
that
while the weighted error initially decreases as $\nu$ decreases, it starts
 increasing when $\nu <2$, justifying the choice $\nu=2$. 
 For the case of $r=16$, similar observations hold
 Although $\nu = 2$  is not the optimal choice when $r=16$, the error for $\nu = 2$
is nearly smallest,
making $\nu = 2$ still a very good candidate for {\bf W-IRKA}.  These numerical results support the idea of choosing $\nu$ and $\varpi$ according
to the decay of the normalized residues. 
Even though this choice seems to yield small weighted errors, there may be variations that are even better. 
 The residues are multiplied by quantities such as $W(-\lambda_i)W(\lambda_i)$, so one might consider incorporating these multiplied quantities as well. 

\begin{figure}[hhh]
\centering
\begin{tabular}{c}
    \epsfxsize=78mm \epsfysize=35mm
\epsffile{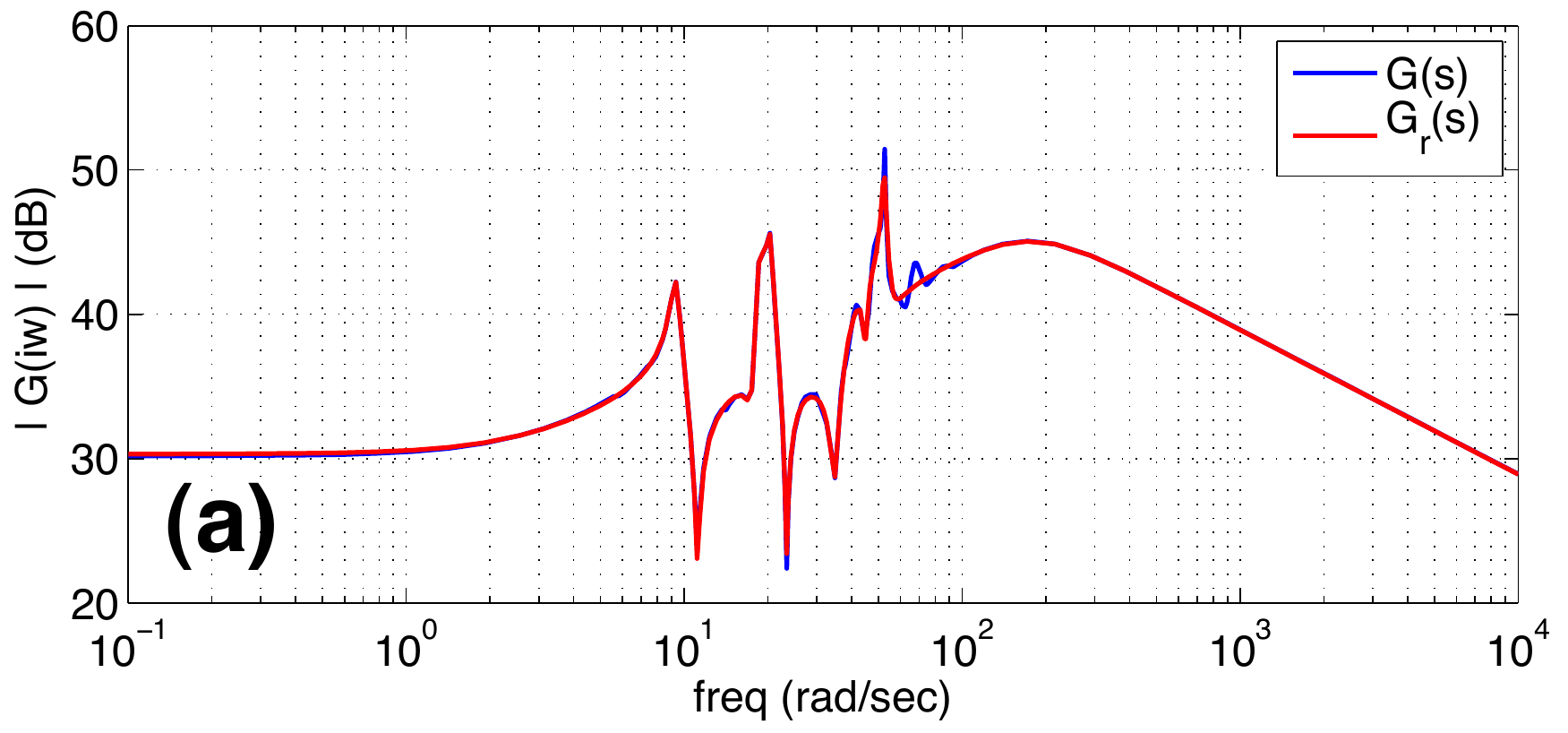}  
\\   \epsfxsize=78mm \epsfysize=35mm
\epsffile{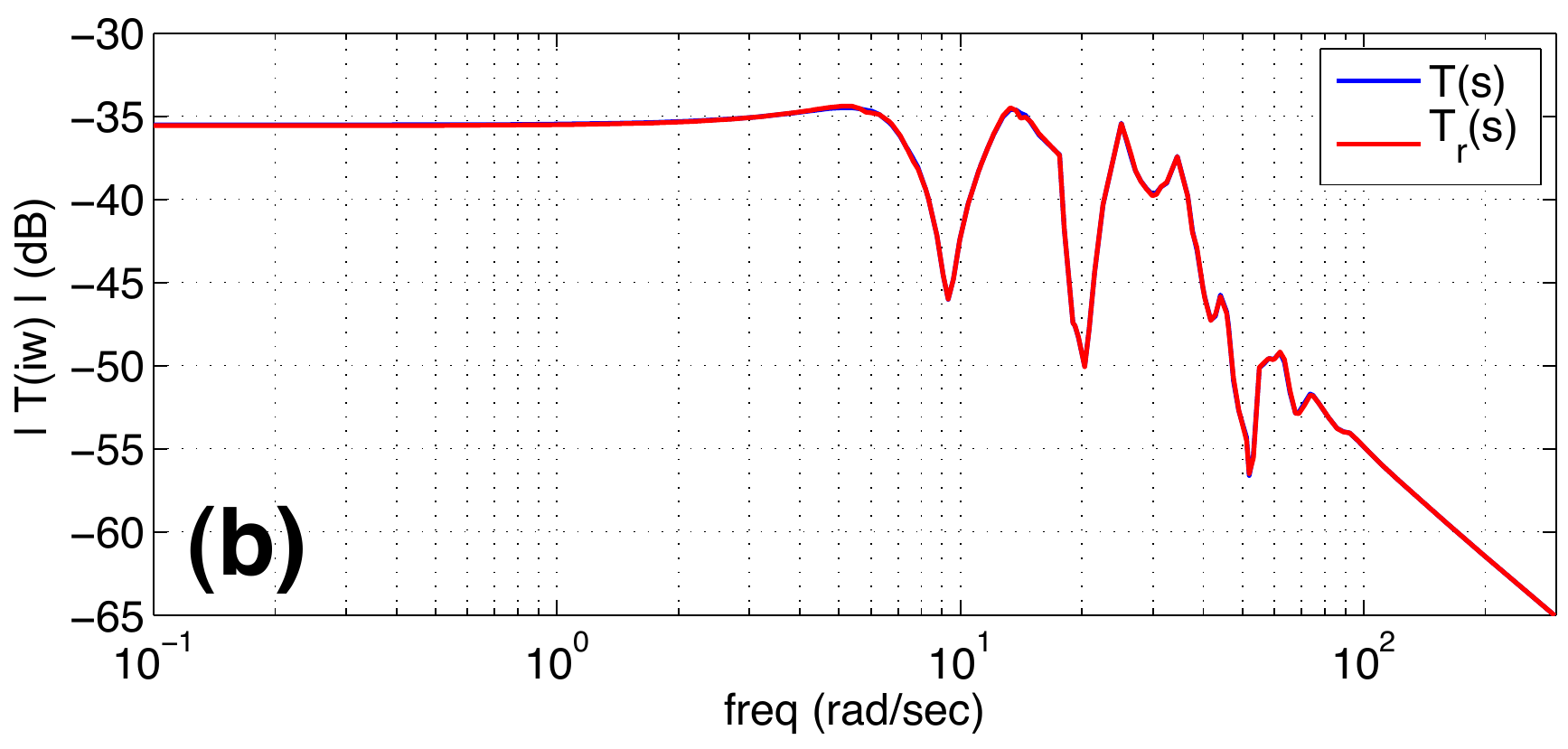}    
\end{tabular}
  \caption{Bode Plots ~(a) Full and reduced controller 
  ~(b) Full and reduced closed-loop system} \label{fig:GandGr}
   \end{figure}

A satisfactory reduced-order controller should not only approximate the full-order controller, but also provide the same closed-loop behavior as the original controller. Let $T$ and $T_r$ denote the full-order and reduced-order closed-loop systems, respectively:  $T$  corresponds to the feedback connection of 
$P$ with $G$; and  $T_r$ to the feedback connection of $P$ with $G_r$. Figure \ref{fig:GandGr}-(a) depicts the amplitude Bode plots  of 
$G$ and $G_r$ for $r=14$ obtained with $\nu = 2$.  $G_r$ is an accurate match to $G$. 
 Figure \ref{fig:GandGr}-(b) shows that the reduced-closed loop behavior $T_r$ almost exactly replicates $T$.
 
 We now compare {\bf W-IRKA} with Frequency Weighted Balanced Truncation ({\bf FWBT}) and
 {\bf IRKA}  of \cite{gugercin2008hmr} for the (unweighted) $\cHtwo$ problem. Comparison with {\bf IRKA} is included to illustrate the importance of including weighting in the $\cHtwo$-based model reduction process.   We vary the reduction order from $r=10$ to $r=20$ in increments of 
$2$, and compute weighted $\cHinf$ and $\cHtwo$ errors for each case. 
We use $\nu=2$ for all cases
even though it might not the best choice for {\bf W-IRKA}. Results are listed in Table \ref{balanced48}. Note that 
for every $r$ value, {\bf W-IRKA} outperforms ${\bf FWBT}$ with respect to the weighted-$\cHtwo$ norm. 
This might be anticipated since {\bf W-IRKA} is designed to reduce the $\cHtwo$ error. 
But {\bf W-IRKA} outperforms ${\bf FWBT}$ with respect to the weighted-$\cHinf$ norm as well in all except 
 the $r=18$ case.  This is significant since balanced truncation approaches generally yield small $\cHinf$ norms.  
 This behavior is similar to the behavior of {\bf IRKA} 
where one often observes that {\bf IRKA}  consistently yields satisfactory $\cHinf$ approximants as well \cite{gugercin2008hmr}. 
Note that 
for $r=10$, the reduced-order controller due to {\bf FWBT}  fails to produce a stable closed-loop system. Table \ref{balanced48} also illustrates that 
 {\bf W-IRKA} significantly outperforms {\bf IRKA} in terms of the weighted error norms. This is what we have expected since unlike {\bf W-IRKA},   {\bf IRKA} is tailored towards the unweighted $\cHtwo$ model reduction problem. 
 This becomes clearer after inspecting Table \ref{IRKAcomp}, which shows that, in terms of the unweighted error $\| G-G_r\|_\cHinf$, {\bf IRKA} outperforms {\bf W-IRKA}. Thus, while $G_r$ 
 from {\bf IRKA} is a better approximation to $G$ in an open-loop sense, once the weight is taken into consideration, {\bf W-IRKA} does what it is designed for, leading to a smaller weighted error.

\begin{table}[hh]
\centering
 $
   ~~~
    {\scriptsize
   \begin{array}{|c|c|c|c|c|c|c|}
   \cline{1-7}
   r & \mbox{10} & \mbox{12} & \mbox{14} &
    \mbox{16} & \mbox{18} & \mbox{20}\\ \cline{1-7}
  \mbox{\bf FWBT}  &1.409 & 0.5286 & 0.0723& 0.0811 &0.0498 & 0.0830\\ \cline{1-7}
  \mbox{\bf W-IRKA}  &0.9175 & 0.1562 & 0.0723&  0.0721 & 0.0722 & 0.0516\\ \cline{1-7}
   \mbox{\bf IRKA} &1.4032 &   0.4837 &   0.1335 &   0.0987 &   0.1194 &   0.1293 \\ \cline{1-7}
  \multicolumn{6}{c}{~~~~~~ \left\|  G-G_r \right\|_{\cHinf(W)}}  \end{array}
  }
$ 
\vspace{2ex}\\
 $
{\scriptsize
  \begin{array}{c|c|c|c|c|c|c|c|c|c|}
    \cline{2-8}
   & \nonumber\mbox{$r$}  & \mbox{10} & \mbox{12} & \mbox{14} &
    \mbox{16} &\mbox{18} & \mbox{20}\\  \cline{2-8}
  & \mbox{\bf FWBT}  & 2.1080 & 1.1723 & 0.1415 &0.1386 &0.1214 & 0.1310\\\cline{2-8}
    &\mbox{\bf W-IRKA} & 0.6677 & 0.2180 & 0.1351 & 0.1309 & 0.1028 &
   0.0956 \\\cline{2-8}
       &\mbox{\bf IRKA} & 1.9540 & 0.8620 & 0.2102 & 0.1286 & 0.2066 &
    0.2317\\\cline{2-8}
 \multicolumn{8}{c}{ \left\|  G-G_r  \right\|_{\cHtwoW}}    \end{array}
 } $
\caption{Comparison of {\bf W-IRKA}, {\bf FWBT}, and {\bf IRKA}}
\label{balanced48}
\end{table}
\vspace{2ex}
\begin{table}[hh]
\centering
 $
   ~~~
    {\scriptsize
   \begin{array}{|c|c|c|c|c|c|c|}
   \cline{1-7}
   r & \mbox{10} & \mbox{12} & \mbox{14} &
    \mbox{16} & \mbox{18} & \mbox{20}\\ \cline{1-7}
  \mbox{\bf W-IRKA} &0.8152  &  0.2750 &   0.3679  &  0.4078 &   0.1274 &   0.0518 \\ \cline{1-7}
   \mbox{\bf IRKA} &0.1062 &   0.1168 &   0.0478 &   0.0513 &   0.0123 &   0.0082 \\ \cline{1-7}
  \multicolumn{6}{c}{ \left\|  G-G_r \right\|_{\cHinf}/ \left\|  G \right\|_{\cHinf} }  \end{array}
  }
$ 
\caption{Comparison of {\bf W-IRKA} and {\bf IRKA}: Unweigthed error}
\label{IRKAcomp}
\end{table}

\subsection{International Space Station 12A Module}
\begin{figure}[htb]
\begin{center}
\subfigure{\includegraphics[width=0.46\textwidth]{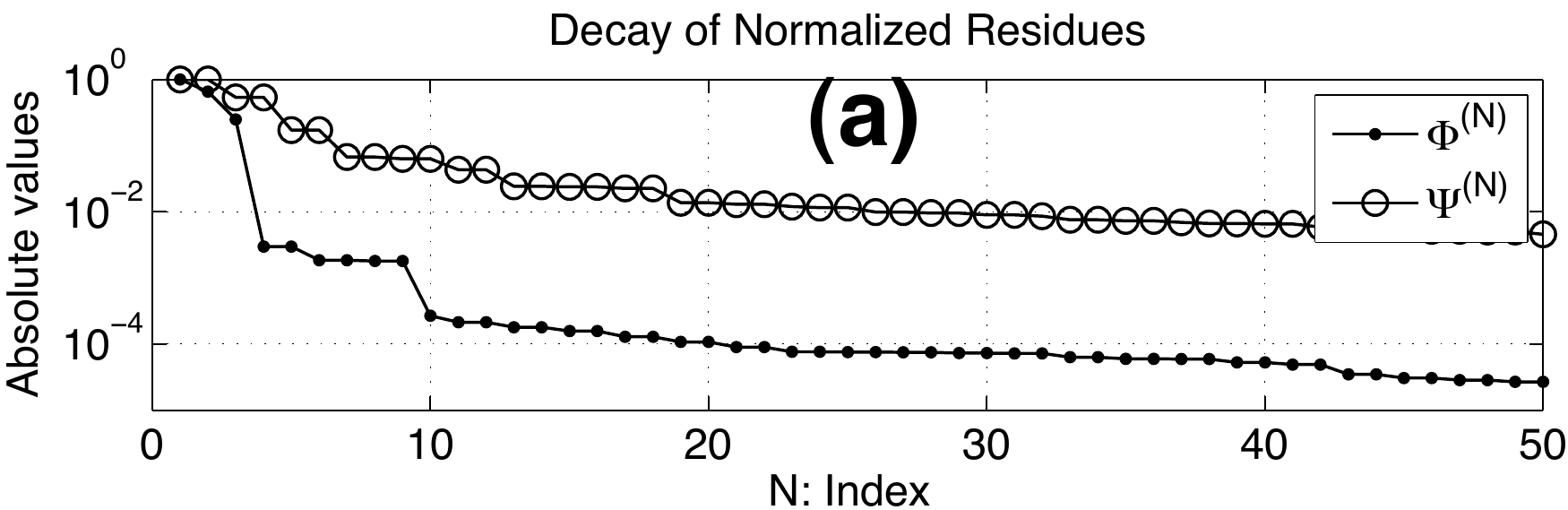}}
\subfigure{\includegraphics[width=0.45\textwidth]{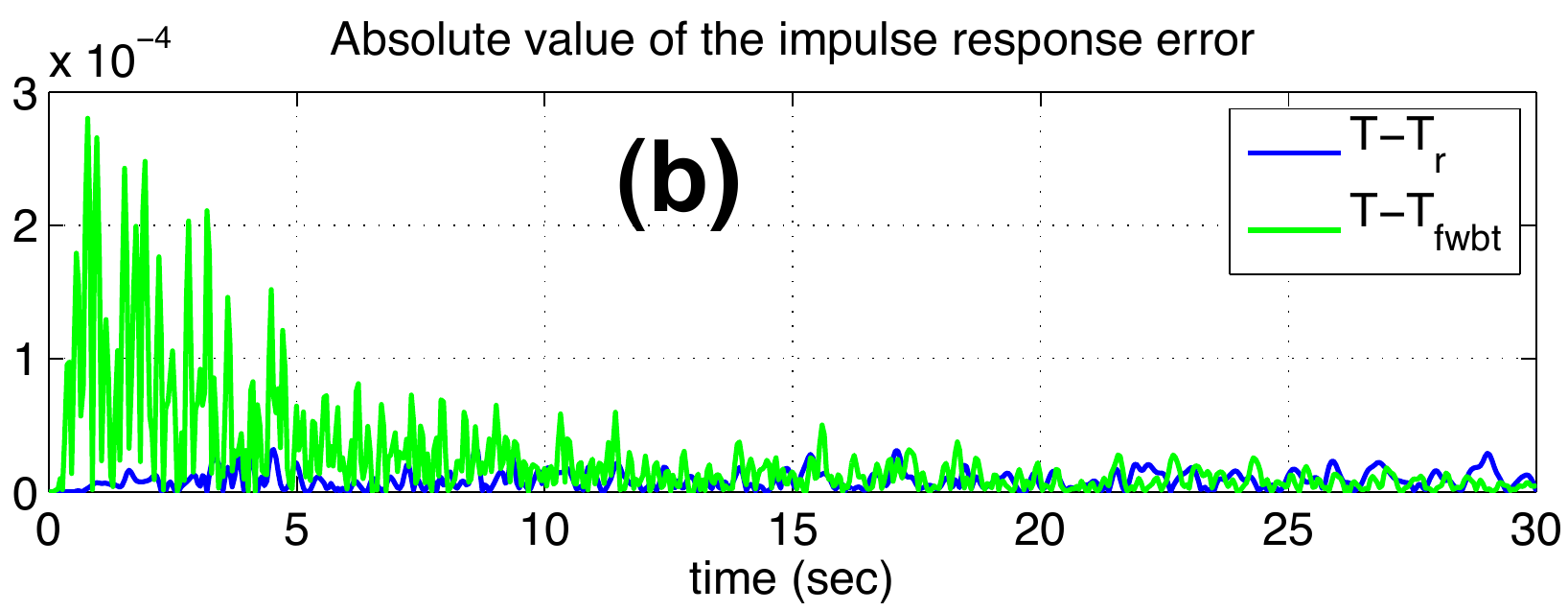}}
\subfigure{\includegraphics[width=0.46\textwidth]{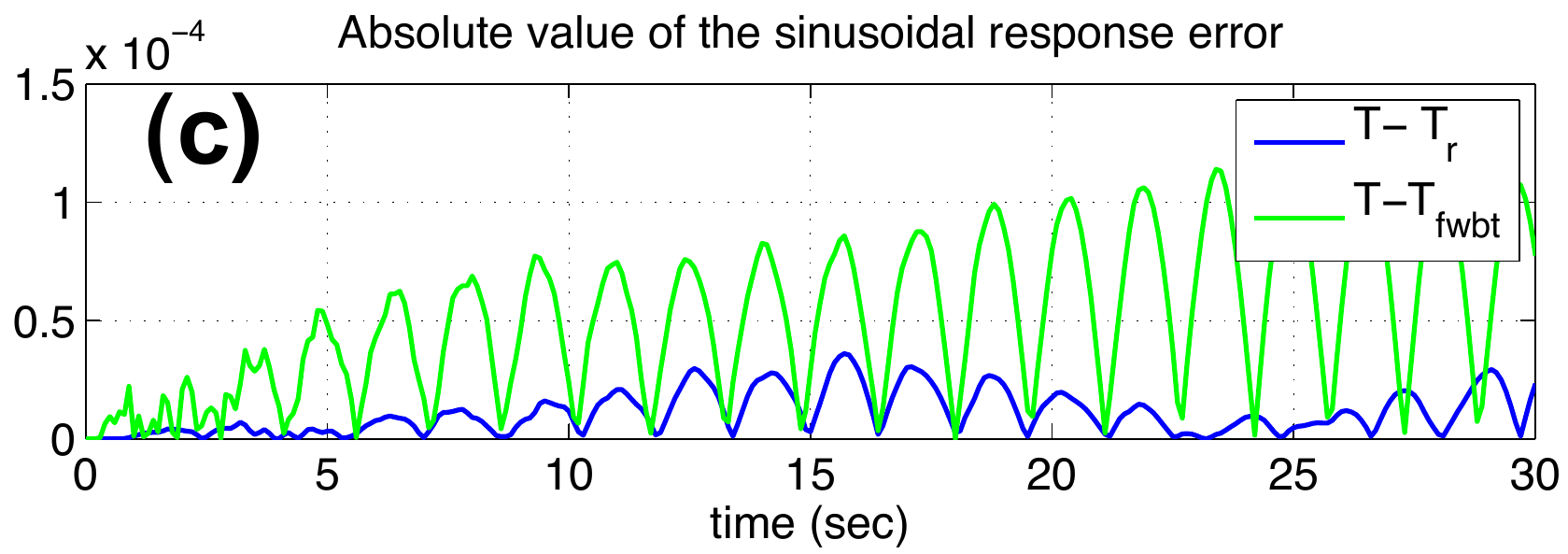}}
\end{center}
\vspace{-2ex} \caption{(a)  Decay of the normalized residues (b)-(c) Comparison of {\bf W-IRKA} and {\bf FWBT} using closed-loop responses.} \label{fig:TandTr}
\end{figure}

The plant, $P$, is a model for the International Space Station 12A Module with dimension $1412$.
It is lightly damped and its impulse response exhibits long-lasting
  oscillations.  A state-feedback,  full-order, observer-based
  controller of order $n = 1412$   is designed to dampen these oscillations. 
  The decay rate of the first $50$ normalized residues $\Phi^{(N)}$ and $\Psi^{(N)}$
 are shown in Figure \ref{fig:TandTr}-(a). While there is almost a two order-of-magnitude drop in $\Phi^{(N)}$ between the third and fourth components,  $\Psi^{(N)}$ continues to stay significant. Hence, we take $\nu=3$ and reduce order from $n=1412$  to $r=60$ using {\bf W-IRKA}. 
For comparison, we also apply {\bf FWBT}.  We denote the resulting reduced-order closed-loop systems due to {\bf W-IRKA} and {\bf FWBT}
by $T_r$ and $T_{\rm fwbt}$, respectively.  Note that $T_{\rm fwbt}$ was unstable for $r=60$.  Indeed, $r=88$ is the \emph{smallest} order {\bf FWBT}-derived reduced controller that lead to a stable closed-system. All {\bf FWBT}-derived $G_r$ are stable; however for $r<88$ when $G_r$ is connected to $P$, the resulting $T_{\rm fwbt}$ is unstable. 
Hence, we compare below the $r=60$ case for {\bf W-IRKA} with the $r=88$ case for {\bf FWBT}.
In Figure \ref{fig:TandTr}-(b), we plot the absolute value of the errors in the impulse responses due to both methods.   {\bf W-IRKA} outperforms {\bf FWBT} even with a lower-order controller.   We also simulate both $T$ and $T_r$ for a sinusoidal input of 
  $u(t) = \cos(2t)$. Results \ref{fig:TandTr}-(c) illustrate the superior performance of {\bf W-IRKA} even more clearly. 

\section{Acknowledgments}     \label{sec:ack}

The work S. Gugercin was supported in part by NSF through Grant DMS-0645347.
The work of A. C. Antoulas was supported in part by NSF through Grant
CCF 1017401, and DFG through Grant 
AN-693/1-1.
 
\section{Conclusions}     \label{sec:conc}

We have presented new formulae for the weighted-$\cHtwo$ inner product and 
 norm that explicitly 
reveal the contribution of  poles and residues both of the full-order model and of the weight.  One of the major consequences of this new representation are new interpolatory optimality conditions for weighted-$\cHtwo$ approximation.  Based on derived weighted-$\cHtwo$ error expressions, we have introduced an approach for producing high-fidelity weighted-$\cHtwo$ reduced models. The effectiveness of this approach has been illustrated with two examples.

\bibliographystyle{plain}        
\bibliography{references}           

\begin{thebibliography}{10}

\bibitem{anderson2002crc}
B.D.O. Anderson and Y.~Liu.
\newblock {Controller reduction: concepts and approaches}.
\newblock {\em IEEE Trans. Automat. Control,}, 34(8):802--812, 2002.

\bibitem{An05}
A.C. Antoulas.
\newblock {\em Approximation of Large-Scale Dynamical Systems}.
\newblock SIAM, 2005.

\bibitem{antoulas2010imr}
A.C. Antoulas, C.A. Beattie, and S.~Gugercin.
\newblock Interpolatory model reduction of large-scale dynamical systems.
\newblock In J.~Mohammadpour and K.~Grigoriadis, editors, {\em Efficient
  Modeling and Control of Large-Scale Systems}. Springer, 2010.

\bibitem{antoulas2001smr}
A.C. Antoulas, D.C. Sorensen, and S.~Gugercin.
\newblock A survey of model reduction methods for large-scale systems.
\newblock {\em Contemporary Mathematics}, 280:193--219, 2001.

\bibitem{devillemagne1987mru}
C.~De~Villemagne and R.E. Skelton.
\newblock {Model reductions using a projection formulation}.
\newblock {\em International Journal of Control}, 46(6):2141--2169, 1987.

\bibitem{enns1984mrw}
D.F. Enns.
\newblock {Model reduction with balanced realizations: An error bound and a
  frequency weighted generalization}.
\newblock In {\em 23rd IEEE Conf. on Decision and Control 1984}, volume~23,
  pages 127--132, 1984.

\bibitem{feldmann1995efficient}
P.~Feldmann and R.W. Freund.
\newblock Efficient linear circuit analysis by {P}ad{\'e} approximation via the
  {L}anczos process.
\newblock {\em IEEE Transactions on Computer-Aided Design of Integrated
  Circuits and Systems}, 14(5):639--649, 1995.

\bibitem{grimme1997krylov}
E.~Grimme.
\newblock {\em {Krylov} Projection Methods for Model Reduction}.
\newblock PhD thesis, Coordinated Science Laboratory, University of Illinois at
  Urbana-Champaign, 1997.

\bibitem{gugercin2008hmr}
S.~Gugercin, A.~C. Antoulas, and C.~A. Beattie.
\newblock $\mathcal{H}_2$ model reduction for large-scale linear dynamical
  systems.
\newblock {\em SIAM J.~Matrix Anal.\ Appl.}, 30(2):609--638, 2008.

\bibitem{gugercin2004smr}
S.~Gugercin and A.C. Antoulas.
\newblock {A survey of model reduction by balanced truncation and some new
  results}.
\newblock {\em International Journal of Control}, 77(8):748--766, 2004.

\bibitem{halevi1992fwm}
Y.~Halevi.
\newblock {Frequency weighted model reduction via optimal projection}.
\newblock {\em IEEE Trans. Automat. Control}, 37(10):1537--1542, 1992.

\bibitem{lin1992mrv}
C.-A. Lin and T.-Y. Chiu.
\newblock {Model reduction via frequency weighted balanced realization}.
\newblock {\em Control Theory and Advanced Technol.}, 8:341--351, 1992.

\bibitem{rommart06}
J.~Rommes and M.~Martins.
\newblock Efficient computation of multivariable transfer function dominant
  poles using subspace acceleration.
\newblock {\em IEEE Trans. on Power Systems}, 21:1471--1483, 2006.

\bibitem{schelfhout2002ncl}
G.~Schelfhout and B.~De~Moor.
\newblock {A note on closed-loop balanced truncation}.
\newblock {\em IEEE Transactions on Automatic Control}, 41(10):1498--1500,
  2002.

\bibitem{spanos1992anewalgorithm}
J.~T. Spanos, M.~H. Milman, and D.~L. Mingori.
\newblock A new algorithm for {$L^2$} optimal model reduction.
\newblock {\em Automatica (Journal of IFAC)}, 28(5):897--909, 1992.

\bibitem{sreeram2005fwm}
V.~Sreeram and A.~Ghafoor.
\newblock Frequency weighted model reduction technique with error bounds.
\newblock In {\em American Control Conference, 2005.}, volume~4, pages 2584 --
  2589, 2005.

\bibitem{varga2003aem}
A.~Varga and B.D.O. Anderson.
\newblock {Accuracy-enhancing methods for balancing-related frequency-weighted
  model and controller reduction}.
\newblock {\em Automatica}, 39(5):919--927, 2003.

\bibitem{wang1999nfw}
G.~Wang, V.~Sreeram, and WQ~Liu.
\newblock {A new frequency-weighted balanced truncation method and an error
  bound}.
\newblock {\em IEEE Trans. on Automat. Control}, 44(9):1734--1737, 1999.

\bibitem{wang2002bpp}
G.~Wang, V.~Sreeram, and WQ~Liu.
\newblock {Balanced performance preserving controller reduction}.
\newblock {\em Systems \& Control Letters}, 46(2):99--110, 2002.

\bibitem{yousouff1984cer}
A.~Yousouff and RE~Skelton.
\newblock {Covariance equivalent realizations with applications to model
  reduction of large-scale systems}.
\newblock {\em Control and Dynamic Systems}, 22:273--348, 1985.

\bibitem{yousuff1985lsa}
A.~Yousuff, D.A. Wagie, and R.E. Skelton.
\newblock {Linear system approximation via covariance equivalent realizations}.
\newblock {\em J. of mathematical analysis and applications}, 106(1):91--115,
  1985.

\bibitem{ZhDG96}
K.~Zhou, J.~Doyle, and K.~Glover.
\newblock {\em Robust and Optimal Control}.
\newblock Prentice-Hall, 1996.

\end{thebibliography}
\end{document}